\documentclass[12pt,reqno]{amsart}
\usepackage{amsmath,amsthm,amsfonts,amssymb,times}
\usepackage{verbatim}
\setbox0=\hbox{$+$}
\newdimen\plusheight
\plusheight=\ht0
\def\+{\;\lower\plusheight\hbox{$+$}\;}

\setbox0=\hbox{$-$}
\newdimen\minusheight
\minusheight=\ht0
\def\-{\;\lower\minusheight\hbox{$-$}\;}

\setbox0=\hbox{$\cdots$}
\newdimen\cdotsheight
\cdotsheight=\plusheight
\def\cds{\lower\cdotsheight\hbox{$\cdots$}}

\renewcommand{\Re}{\operatorname{Re}}

\makeatletter
\def\leqalignno#1{\displ@y \tabskip\z@ plus\@ne fil
  \halign to\displaywidth{\hfil$\@lign\displaystyle{##}$\tabskip\z@skip
    &$\@lign\displaystyle{{}##}$\hfil\tabskip\z@ plus\@ne fil
    &\kern-\displaywidth\rlap{$\@lign\hbox{\rm##}$}\tabskip\displaywidth\crcr
    #1\crcr}}
\makeatother

\newcommand{\eb}{\begin{equation}}
\newcommand{\ee}{\end{equation}}

\newcommand{\df}{\dfrac}
\newcommand{\tf}{\tfrac}

\renewcommand{\Re}{\operatorname{Re}}

\newcommand{\G}{\Gamma}

\renewcommand{\Re}{\text{Re}}

\renewcommand{\(}{\left\(}
\renewcommand{\)}{\right\)}
\renewcommand{\[}{\left\[}
\renewcommand{\]}{\right\]}
\renewcommand{\i}{\infty}

\numberwithin{equation}{section}
 \theoremstyle{plain}

\newtheorem{theorem}{Theorem}[section]
\numberwithin{equation}{section}
\theoremstyle{plain}
\newtheorem{lemma}[theorem]{Lemma}

\newtheorem{entry}[theorem]{Entry}
\newtheorem{question}[theorem]{Question}

\allowdisplaybreaks

\begin{document}

\title[Beautiful Integrals]
{Ramanujan's Beautiful Integrals}
\author{Bruce C.~Berndt, Atul Dixit}
\address{Department of Mathematics, University of Illinois, 1409 W.~Green St., Urbana, IL 61801, USA}
\address{Department of Mathematics, Indian Institute of Technology Gandhinagar, Palaj, Gandhinagar 382355, Gujarat, India}
\maketitle

\begin{abstract} Throughout his entire mathematical life, Ramanujan loved to evaluate definite integrals.  One can find them in his problems submitted to the \emph{Journal of the Indian Mathematical Society}, notebooks, Quarterly Reports to the University of Madras,  letters to Hardy, published papers and the Lost Notebook. His evaluations are often surprising, beautiful, elegant, and useful in other mathematical contexts. He also discovered general methods for evaluating and approximating integrals.  A survey of Ramanujan's contributions to the evaluation of integrals is given, with examples provided  from each of the above-mentioned sources.

\end{abstract}

\section{Introduction}
Ramanujan loved infinite series and integrals.  They permeate almost all of his work from the years he recorded his findings in notebooks \cite{nb} until the end of his life in 1920 at the age of 32. In this paper we provide a survey of some of his most beautiful theorems on integrals. Of course, it is impossible to adequately cover what Ramanujan accomplished in his devotion to integrals.   Many of Ramanujan's theorems and examples of integrals have inspired countless mathematicians to take Ramanujan's thoughts and proceed further. For many of Ramanujan's integrals, we stand in awe and admire their beauty, much as we listen to a beautiful Beethoven piano sonata or an intricate but mellifluous raaga in Carnatic or Hindustani classical music.  We hope that this survey will provide further inspiration.

 Ramanujan evaluated many definite integrals, most often infinite integrals.  In many cases, the integrals are so ``unusual,'' that we often wonder how Ramanujan ever thought that elegant evaluations existed.  Some of his integrals satisfy often surprising functional equations.  He was an expert in finding exquisite examples for integral transforms, some of which are original with him.  His so-called ``Master Theorem" fits into this category.  Some of his integrals have (non-trivial) relations with infinite series and continued fractions.  Ramanujan was also a master in finding asymptotic expansions of integrals.

 Integrals often arise in concomitant problems that Ramanujan studied.  For example, in his first letter to G.~H.~Hardy, Ramanujan asserted \cite[p.~24]{bcbrar} (with a minor correction needed)
 \begin{quote}
1, 2, 4, 5, 8, 9, 10, 13, 16, 17, 18,\dots\, are numbers which are either themselves squares or which can be expressed as the sum of two squares.
The number of such numbers greater than \emph{A}  and less than
\emph{B}
 $$=K\int_A^B \df{dx}{\sqrt{\log x}} + \theta (x)$$
where $ K  = .764 \dots$ and $\theta (x)$ is very small when
compared with the previous integral.  \emph{K} and $\theta (x)$
have been exactly found though complicated.
\end{quote}

Theorems and claims of this kind are better addressed in the contexts in which they arise, and so we do not address such integral appearances in the present paper.

Beginning in 1911, Ramanujan offered a total of 58 \emph{Questions} to the \emph{Journal of the Indian Mathematical Society}.  In seven of them, readers are asked to evaluate definite integrals.  (For discussions of all 58 problems, see a paper by the first author, Y.-S.~Choi, and S.-Y.~Kang \cite{lange}, \cite[pp.~215--258]{bcbrara}.) Ramanujan's first letter to Hardy contains over 60 statements, and eleven of them pertain to definite integrals. In his second letter to Hardy, seven entries provide evaluations of integrals \cite[pp.~xxiii--xxix]{cp}, \cite[Chapter 2]{bcbrar}.  (Portions of both letters have been lost.)  Six of Ramanujan's published papers are devoted to the evaluation of integrals.  The most abundant source for Ramanujan's integrals are his (earlier) notebooks \cite{nb}.  His lost notebook \cite{lnb} also contains several intriguing integrals \cite{lnb}.

The purpose of this paper is to provide readers with a survey of some of Ramanujan's most beautiful, surprising, and possibly useful integrals.  Examples from each of the sources mentioned in the previous paragraph will be given.

\section{Problems Posed by Ramanujan in the\\ \emph{Journal of the Indian Mathematical Society}}

Question 783 below is an especially elegant problem posed by Ramanujan in the \emph{Journal of the Indian Mathematical Society} \cite{783} and found in Ramanujan's third notebook \cite[p.~373]{nb}.  Like so many of Ramanujan's discoveries, we wonder how Ramanujan ever thought of this problem.  While reading this paper, readers will undoubtedly ask questions of this sort several times.

\begin{question} \textbf{(Question 783)}\label{783a}
For $n\geq 0$, put $v=u^n-u^{n-1}$, and define
\begin{equation}\label{783b}
\varphi(n):=\int_0^1\df{\log u}{v}dv.
\end{equation}
Then,
\begin{equation*}
\varphi(0)=\df{\pi^2}{6},\qquad \varphi(1)=\df{\pi^2}{12}, \qquad\text{and}\qquad\varphi(2)=\df{\pi^2}{15}.
\end{equation*}
Furthermore, if $n>0$,
\begin{equation}\label{783c}
\varphi(n)+\varphi\left(\frac{1}{n}\right)=\df{\pi^2}{6}.
\end{equation}
\end{question}

Inspired by Question \ref{783a}, Berndt and R.~J.~Evans \cite{rje} established the generalization given below.  Here and in the sequel, we use the conventions
\begin{equation*}
f(\infty):=\begin{cases}\underset{x\to\infty}{\lim}f(x), \quad &\text{provided that the limit exists},\\
 \infty, \quad & \text{if } f(x)\to \infty \text{ as } x\to\infty.
 \end{cases}
 \end{equation*}

\begin{theorem}\label{783d}
Let $g$ be a strictly increasing, differentiable function on $[0,\infty)$ with $g(0)=1$ and $g(\infty)=\infty$. For $n>0$ and $t\geq0$, define
\begin{equation*}
v(t):=\df{g^n(t)}{g(1/t)}.
\end{equation*}
Suppose that
\begin{equation*}
\varphi(n):=\int_0^1\log\,g(t)\df{dv}{v}
\end{equation*}
converges. Then
\begin{equation*}
\varphi(n)+\varphi\left(\frac{1}{n}\right)=2\varphi(1).
\end{equation*}
\end{theorem}

Note that if $g(t)=1+t$, Theorem \ref{783d} reduces to Question \ref{783a}. See \cite{rje} and \cite[pp.~326--329]{IV} for more details.

The integral \eqref{783b} is reminiscent of the dilogarithm, defined by
\begin{equation}\label{dilog}
\text{Li}_2(z):=-\int_0^z\df{\log(1-w)}{w}dw, \qquad z\in\mathbb{C},
\end{equation}
where the principal branch of $\log(1-w)$ is chosen.  The dilogarithm was studied by Ramanujan in Chapter 9 of his first notebook \cite[Entries 5--7]{nb}, \cite[pp.~246--249]{I}, where many of the fundamental properties of $\text{Li}_2(z)$ are proved.  He returned to $\text{Li}_2(z)$ in his third notebook \cite[pp.~365]{nb},  \cite[pp.~322--326]{IV}.

Ramanujan was fond of formulas evincing symmetry, such as in Question 295 below \cite{295}.

\begin{question}\label{295a} If $\alpha$ and $\beta$ are positive and $\alpha\beta=\pi$, then
\begin{equation}\label{1}
\sqrt{\alpha}\int_0^{\i}\df{e^{-x^2}}{\cosh \, \alpha x}dx =\sqrt{\beta}\int_0^{\i}\df{e^{-x^2}}{\cosh \, \beta x}dx.
\end{equation}
\end{question}

 Question \ref{295a} can be found in Ramanujan's first letter to G.~H.~Hardy \cite[p.~27]{bcbrar}, and also in Section 21 of Chapter 13 in his second notebook \cite{nb}, \cite[p.~225]{II}.  The identity \eqref{1} appears in a manuscript published with the lost notebook \cite{lnb}, \cite[p.~368]{geabcbIV}, and it  was also  established by Hardy \cite[p.~203]{hardy1}.
  Four additional relations in the spirit of \eqref{1} can be found in the aforementioned manuscript \cite[p.~368]{geabcbIV}.

  \section{Ramanujan's First Two Letters to Hardy} In his first letter to Hardy, Ramanujan offers the following ``reciprocity" formula or ``theta" relation.

  \begin{theorem}\label{modular1} For $n>0$, define
  \begin{equation}\label{modular}
 \phi(n):=\int_0^{\infty} \dfrac{\cos nx}{ e^{2\pi \sqrt x} - 1}dx.
\end{equation}
 Then
 \begin{equation}\label{modular2}
 \int_0^{\infty} \df{\sin nx}{e^{2\pi \sqrt
x} - 1}dx = \phi (n) - \df{1}{2n} + \phi \left(\df{\pi ^2}{
n}\right)\sqrt{\df{2 \pi ^3}{ n^3}} .
\end{equation}
``$\phi (n)$ is a complicated function \dots"
\end{theorem}

As special cases,

$$\phi (0) = \df{1}{12};\qquad \phi \left(\df{\pi}{ 2}\right) = \df{1}{4\pi};\qquad \phi (\pi) = \df{2 -\sqrt{2}}{ 8};\qquad \phi(2\pi) = \df{1}
{16};$$

$$\phi \left(\df{2\pi}{ 5}\right) = \df{8 - 3\sqrt{5}}{16};\qquad \phi
\left(\df{\pi}{ 5}\right) = \df{6 +\sqrt{5}}{ 4} - \df{ 5 \sqrt{10}}{8};\qquad \phi
(\infty) = 0;$$

$$\phi \left(\df{2\pi }{ 3}\right) = \df{1}{ 3} - \sqrt{3} \left(\df {3}{ 16} - \df{1}{
8 \pi}\right) .$$

An abbreviated version of Theorem \ref{modular1} appeared as a problem in the \emph{Journal of the Indian Mathematical Society} \cite{463}.  Theorem \ref{modular1} can also be found in Ramanujan's notebooks \cite{nb}.

The function $\phi(n)$  in \eqref{modular} can be expressed in terms of a variant of Gauss sums to which Ramanujan devotes an entire paper \cite{12}.  For example, if $n=a/b$, where $a$ and $b$ are positive odd numbers, then
\begin{equation}\label{modular3}
\phi\left(\df{\pi a}{b}\right):=\df14\sum_{r=1}^b(b-2r)\cos\left(\df{r^2\pi a}{b}\right)-\df{b}{4a}\sqrt{\df{b}{a}}
\sum_{r=1}^{a}(a-2r)\sin\left(\df{1}{4}\pi+\df{r^2 \pi b}{a}\right).
\end{equation}
See Berndt's book \cite[pp.~296--303]{IV} for a more complete discussion of $\phi(n)$ and its  connections with Ramanujan's analogues of Gauss sums and another class of his infinite series.

  After discussing the Rogers--Ramanujan continued fraction and two generalizations in his second letter to Hardy, Ramanujan offers representations for a pair of integrals by continued fractions \cite[p.~xxviii]{cp}, \cite[p.~57]{bcbrar}.

  \begin{theorem} We have
  \begin{align*}
  4\int_0^{\infty}\df{xe^{-x\sqrt5}}{\cosh x}dx&=\df{1}{1}\+\df{1^2}{1}\+\df{1^2}{1}\+\df{2^2}{1}\+\df{2^2}{1}\+\df{3^2}{1}\+\df{3^2}{1}\+\cds,\\
2\int_0^{\infty}\df{x^2e^{-x\sqrt3}}{\sinh x}dx&=\df{1}{1}\+\df{1^3}{1}\+\df{1^3}{3}\+\df{2^3}{1}\+\df{2^3}{5}\+\df{3^3}{1}\+\df{3^3}{7}\+\cds.
\end{align*}
\end{theorem}

These two identities were first proved in print by C.~T.~Preece \cite{preece} in 1931.  Above, we have corrected two misprints in the second formula that appears in \cite[p.~xxviii]{cp}.

\section{Ramanujan's Published Papers on Integrals}

Ramanujan published the papers \cite{7}, \cite{11}, \cite{12}, \cite{22}, \cite{23}, and \cite{27} on definite integrals.  We briefly discussed some of the content of \cite{12} in the previous section.

Define, for $\Re(w)\geq 0$,
\begin{equation*}
\phi_{w}(t):=\int_0^{\i}\df{\cos\,\pi tx}{\cosh\,\pi x}e^{-\pi wx^2}dx \qquad\text{and}\qquad
\psi_{w}(t):=\int_0^{\i}\df{\sin\,\pi tx}{\sinh\,\pi x}e^{-\pi wx^2}dx.
\end{equation*}
 In \cite{23}, Ramanujan made use of ``modular relations'' satisfied by $\phi_w(t)$ and $\psi_w(t)$, namely,
 \begin{align*}
   \phi_w(t)&=\df{1}{\sqrt{w}}e^{-\tf14\pi t^2/w}\psi_{1/w}(it/w)
   \end{align*}
   and
   \begin{equation*}
   e^{\tf14\pi t^2/w}\left\{\frac12+\psi_w(t)\right\}=e^{\tf14\pi(t+w)^2/w}\phi_w(t+w),
   \end{equation*}
   to develop representations in terms of theta functions.  These were then used
     to evaluate  large classes of integrals for specific values of $t$ and $w$.  We give two examples:
\begin{align}\label{mos}
\int_0^{\i}\df{\sin\,2\pi tx}{\sinh\,\pi x}\cos\,\pi x^2\,dx&=\df{\cosh\,\pi t-\cos\,\pi t^2}{2\sinh\,\pi t},\\\
\int_0^{\i}\df{\sin\,2\pi tx}{\sinh\,\pi x}\sin\,\pi x^2\,dx&=\df{\sin\,\pi t^2}{2\sinh\,\pi t}.\notag
\end{align}
The integral evaluation \eqref{mos} was used by A.~K.~Mustafy \cite{mustafy} to obtain a new integral representation of the Riemann zeta function $\zeta(s)$, with an integrand similar to those above, which in turn also gives a proof of Riemann's functional equation for $\zeta(s)$.

Perhaps Ramanujan's most important paper on integrals is \cite{27}. Here, several classes of integrals involving the gamma function are evaluated in closed form.   Many of the integrals in \cite{27} can perhaps be evaluated by contour integration, although Ramanujan did not use this method in \cite{27}.
It  has generally been accepted that Ramanujan was not conversant with the analytic theory of functions of a complex variable. Hardy opined
\cite[p.~xxx]{cp}, ``\dots and he (Ramanujan) had indeed but the vaguest idea of what a function of a complex variable was.'' However, in \cite{27} it is clear that Ramanujan knew certain elementary facts about functions of a complex variable, in particular, at the very least, he was aware of the care needed in treating branches of ``multi-valued'' functions.
 On the next-to-last page of Ramanujan's third notebook \cite[pp.~391]{nb}, several integrals from complex analysis are recorded.  Next to one of them appear the words, ``contour integration.''   So, maybe Ramanujan knew more complex analysis than either Hardy or others have thought.

One of the general integrals that Ramanujan evaluated is \cite{27}, \cite[pp.~221--222]{cp}:
\begin{equation*}
\int_{-\i}^{\i}\Gamma(\alpha+x)\Gamma(\beta-x)e^{inx}dx.
\end{equation*}
Similarly, Ramanujan devised a general approach to evaluating
\begin{equation*}
\int_{-\i}^{\i}\df{e^{inx}}{\Gamma(\alpha+x)\Gamma(\beta-x)\Gamma(\gamma+\ell x)\Gamma(\delta-\ell x)}dx,
\end{equation*}
where $n$ and $\ell$ are real numbers.
We forego hypotheses in this brief survey, but instead mention only a special case.  If $\alpha+\beta+\gamma+\delta=4$, then \cite[p.~229]{27}
\begin{gather*}
\int_{-\i}^{\i}\df{\cos\{\pi(x+\beta+\gamma)\}}{\Gamma(\alpha+x)\Gamma(\beta-x)\Gamma(\gamma+2 x)\Gamma(\delta-2x)}dx\\=
\df{1}{2\Gamma(\gamma+\delta-1)\Gamma(2\alpha+\delta-2)\Gamma(2\beta+\gamma-2)}.
\end{gather*}

Ramanujan  devised some beautiful and unusual definite integral evaluations involving products of ordinary Bessel functions $J_{\nu}(x)$.  For example, if
$\Re(\alpha+\beta)>-1$ \cite[p.~225]{27},
\begin{equation*}
\int_{-\i}^{\i}\df{J_{\alpha+w}(x)}{x^{\alpha+w}}\df{J_{\beta-w}(y)}{y^{\beta-w}}dw=
\df{J_{\alpha+\beta}\{\sqrt{(2x^2+2y^2)}\}}{(\tf12 x^2+\tf12 y^2)^{(\alpha+\beta)/2}}.
\end{equation*}
Ramanujan concluded his paper \cite{27} with a formula providing the evaluation of a fairly general integral involving the product of four Bessel functions.

Ramanujan's integrals were discussed by Watson in his treatise \cite[p.~449]{watson}.  In a footnote, he remarks, ``these integrals evaluated by Ramanujan may prove to be of the highest importance in the theory of the transmission of Electric Waves.''

One of the highly influential papers that Ramanujan wrote after arriving in England is his paper \cite{riemann}, in which he obtains transformations for integrals arising in the theory of Riemann's zeta function $\zeta(s)$.  The Riemann $\xi(s)$ and $\Xi$ functions are defined by \cite[p.~16]{titch}
\begin{align*}
\xi(s):=\frac{1}{2}s(s-1)\pi^{-\frac{s}{2}}\Gamma\left(\frac{s}{2}\right)\zeta(s) \qquad\text{and}\qquad \Xi(t)&:=\xi\left(\tfrac{1}{2}+it\right).
\end{align*}
In his review \cite{hardyw} of Ramanujan's work in England until 1917, Hardy cites \cite{riemann} as one of Ramanujan's four most important papers.

One of the two main results in his paper is \cite[Equation (12)]{riemann}
\begin{align}\label{uncanny1}
&\int_{0}^{\infty}\left\{e^{-z}-4\pi\int_{0}^{\infty}\frac{xe^{-3z-\pi x^2e^{-4z}}}{e^{2\pi x}-1}\, dx\right\}\cos(tz)\, dz\nonumber\\
&=\frac{1}{8\sqrt{\pi}}\G\left(\frac{-1+it}{4}\right)\G\left(\frac{-1-it}{4}\right)\Xi\left(\frac{t}{2}\right),
\end{align}
 which, through Fourier inversion, leads to \cite[Equation (13)]{riemann}
\begin{align}\label{13}
&e^{-n}-4\pi e^{-3n}\int_{0}^{\infty}\frac{xe^{-\pi x^2e^{-4n}}}{e^{2\pi x}-1}\, dx\nonumber\\
&=\frac{1}{4\pi\sqrt{\pi}}\int_{0}^{\infty}\G\left(\frac{-1+it}{4}\right)\G\left(\frac{-1-it}{4}\right)\Xi\left(\frac{t}{2}\right)\cos(nt)\, dt
\end{align}
for $n\in\mathbb{R}$. About \eqref{13}, Hardy \cite{ghh} writes, for $\sigma=\Re(s)$,

\begin{quote}\emph{The integral has properties similar to those of the integral by means of which I proved recently that $\zeta(s)$ has an infinity of zeros on the line $\sigma=1/2$, and may be used for the same purpose.}\end{quote}

  In the last section of his paper \cite{riemann}, for  $n\in\mathbb{R}$, Ramanujan obtains new beautiful integral representations for
\begin{equation}\label{sec5int}
F(n,s):=\int_{0}^{\infty}\Gamma\left(\frac{s-1+it}{4}\right)\Gamma\left(\frac{s-1-it}{4}\right)
\Xi\left(\frac{t+is}{2}\right)\Xi\left(\frac{t-is}{2}\right)\frac{\cos nt}{(s+1)^2+t^2}\, dt,
\end{equation}
 each valid in some vertical strip in the half-plane $\sigma>1$. One of these is given by
\begin{align*}
F(n, s)=\frac{1}{8}(4\pi)^{-\frac{1}{2}(s-3)}\int_{0}^{\infty}x^{s}\left(\frac{1}{\exp{(xe^n)}-1}
-\frac{1}{xe^n}\right)\left(\frac{1}{\exp{(xe^{-n})}-1}-\frac{1}{xe^{-n}}\right)\, dx.
\end{align*}
(See \cite{dixit} for corrected misprints in \cite{riemann}.)
Regarding the special case $s=0$ of this integral, Hardy \cite{ghh} writes,
 \begin{quote}\emph{\dots the properties of this integral resemble those of one which Mr. Littlewood and I have used, in a paper to be published shortly in Acta Mathematica to prove that}
 \end{quote}
\begin{equation*}
\int_{-T}^{T}\left|\zeta\left(\frac{1}{2}+ti\right)\right|^2\, dt \sim
2 T\log T\hspace{3mm}(T\to\infty),
\end{equation*}\\[2pt]
 where a misprint from Hardy's paper \cite{ghh} has been corrected.
This special case $s=0$ of the integral in \eqref{sec5int} also appears on page $220$ of the lost notebook \cite{lnb}, where Ramanujan gives an exquisitely beautiful modular relation associated with it. The reader is referred to \cite{bcbad}, \cite{dixitms} and \cite{dixzah} for more details on Ramanujan's formulas from \cite{riemann}, their importance and their applications.

In conclusion, about Ramanujan's formulas from \cite{riemann}, Hardy remarks \cite{ghh},
\begin{quote}\emph{It is difficult at present to estimate the importance of these results. The unsolved problems concerning the zeros of $\zeta(s)$ or of $\Xi(t)$ are among the most obscure and difficult in the whole range of Pure Mathematics. \dots
  But I should not be at all surprised if still more important applications were to be made of Mr. Ramanujan's formulae in the future.}
\end{quote}

\section{Ramanujan's Quarterly Reports}

 Ramanujan's fame began with his publications in the \emph{Journal of the Indian Mathematical Society} in 1911,
 and it reached the English astronomer Sir Gilbert Walker, who was working at an observatory in Madras.  In a letter to the University of Madras dated February 26, 1913, he wrote, ``The University would be justified in enabling S.~Ramanujan for a few years \emph{at least} to spend the whole of his time on mathematics, without any anxiety as to his livelihood.'' The Board of Studies at the University of Madras agreed to this request, and its chairman, Professor B.~Hanumantha Rao, wrote a letter to Vice-Chancellor Francis Dewsbury on March 25, 1913 with the recommendation that Ramanujan be awarded a scholarship of 75 rupees per month.   A stipulation in the scholarship required Ramanujan to write Quarterly Reports to the Board of Studies in Mathematics. Ramanujan wrote three of these Quarterly Reports before he departed for England on March 17, 1914.  Unfortunately, they were eventually lost; but, on the other hand,  fortunately, T.~A.~Satagopan made a handwritten copy of the reports in 1925.  An extensive description of their contents was published by Berndt in the \emph{Bulletin of the London Mathematical Society} \cite{quartblms}. The aforementioned letters can be found in the book \cite[pp.~70--76]{bcbrar} by Berndt and Robert Rankin.

 We offer a few results from the Quarterly Results; the first are Frullani's Theorem and Ramanujan's generalization.

\begin{theorem}\label{frullani1} \textbf{(Frullani)} If $f$ is a continuous function on $[0,\infty)$ such that $f(\infty)$ exists, then, for any pair $a,b>0$,
\begin{equation}\label{frullani}
\int_0^{\infty}\df{f(ax)-f(bx)}{x}dx=\{f(0)-f(\infty)\}\log\left(\df{b}{a}\right).
\end{equation}
If $f(\infty)$ does not exist, but $f(x)/x$ is integrable over $[c,\infty)$ for  $c>0$, then \eqref{frullani} still holds, but with $f(\infty)$ replaced by 0.
\end{theorem}

In his second Quarterly Report, Ramanujan offers a beautiful generalization of Frullani's Theorem. A slightly less general version is provided by Ramanujan in the unorganized pages of his second notebook \cite[pp.~332, 334]{nb}, \cite[p.~316]{I}. We do not give below the hypotheses that are needed for $u(x)$ and $v(x)$;  see \cite[pp.~299, 313]{I} for these requirements.  Set
$$f(x)-f(\infty)=\sum_{k=0}^{\infty}\df{u(k)(-x)^k}{k!}\qquad\text{and}\qquad g(x)-g(\infty)=\sum_{k=0}^{\infty}\df{v(k)(-x)^k}{k!}.$$
 Ramanujan also assumes that the limit below can be taken under the integral sign.

 \begin{theorem}\label{frullani2}
 Let  $u(x)$ and $v(x)$ be given as above, and assume that $f$ and $g$ are continuous functions on $[0,\infty)$.  Also assume that $f(0)=g(0)$ and $f(\infty)=g(\infty)$.  Then, if $a,b>0$,
 \begin{equation*}
 \lim_{n\to0}\int_0^{\infty}x^{n-1}\left\{f(ax)-g(bx)\right\}dx
 =\{f(0)-f(\infty)\}\left\{\log\left(\df{b}{a}\right)+\df{d}{ds}\left(\log\left(\df{v(s)}{u(s)}\right)\right)_{s=0}\right\}.
 \end{equation*}
 \end{theorem}

 Ramanuan's proof depends on his now famous \emph{Master Theorem}.  He  assumes that a function $F(x)$ can be expanded in a Taylor series about $x=0$ with an infinite radius of convergence.  Then Ramanujan asserts that the value of the integral
$$ \int_0^{\infty}x^{n-1}F(x)dx$$
can be found from the coefficient of $x^n$ in the expansion of $F(x)$.

\begin{theorem}\label{mastertheorem} \textbf{Ramanujan's Master Theorem.} Suppose that for $-\infty<x<\infty$,
\begin{equation}\label{master}
F(x)=\sum_{k=0}^{\infty}\df{\varphi(k)(-x)^k}{k!}.
\end{equation}
Then
\begin{equation}\label{master1}
\int_0^{\infty}x^{n-1}\sum_{k=0}^{\infty}\df{\varphi(k)(-x)^k}{k!}dx=\Gamma(n)\varphi(-n).
\end{equation}
\end{theorem}

For our first illustration of Ramanujan's Master Theorem \cite[p.~300]{I}, let $m,n>0$ and set $x=y/(1+y)$ to find that
\begin{equation*}\int_0^1x^{m-1}(1-x)^{n-1}dx=\int_0^{\infty}y^{m-1}(1+y)^{-(m+n)}dy.\end{equation*}
From the binomial series
$$ (1+y)^{-r}=\sum_{k=0}^{\infty}\df{\Gamma(k+r)}{\Gamma(r)k!}(-y)^k, \qquad |y|<1,$$
we find that $\varphi(t)=\Gamma(t+m+n)/\Gamma(m+n).$  Applying Ramanujan's Master Theorem, we deduce the well-known representation of the beta function $B(m,n)$,
\begin{equation}\label{beta}
B(m,n):=\int_0^1 x^{m-1}(1-x)^{n-1}dx=\Gamma(m)\varphi(-m)=\df{\Gamma(m)\Gamma(n)}{\Gamma(m+n)}.
\end{equation}

For our second example, we need the notation
\begin{equation*}\label{intro1}
(a;q)_0:=1,\quad(a;q)_n:=\prod_{k=0}^{n-1}(1-aq^k),\quad n\geq1,
 \end{equation*}
 and
 \begin{equation*}\label{intro2}
 (a;q)_{\i}:=\lim_{n\to\i}(a;q)_n, \quad |q|<1.
 \end{equation*}

Recall the $q$-binomial theorem \cite[p.~8]{gasper}
 \begin{equation*}
 \sum_{m=0}^{\i}\df{(a;q)_m}{(q;q)_m}z^m=\df{(az;q)_{\i}}{(z;q)_{\i}},\qquad |z|<1.
 \end{equation*}
 Letting $z=-x$, replacing $a$ by $aq$, and applying the Master Theorem, we establish,  for $\Re\,s>0$, Ramanujan's beautiful identity \cite{12}, \cite[p.~57]{cp},
 \begin{equation}\label{beautiful}
 \int_0^{\infty}t^{s-1}\df{(-atq;q)_{\infty}}{(-t;q)_{\infty}}dt=\df{\pi}{\sin(\pi s)}\df{(q^{1-s};q)_{\infty}(aq;q)_{\infty}}{(q;q)_{\infty}(aq^{1-s};q)_{\infty}}.
 \end{equation}
 Richard Askey \cite{askey} made a thorough study of the integral in \eqref{beautiful} and showed that, if $s=x$ and $a=q^{x+y}$, then this integral is a natural $q$-analogue of the beta function $B(x,y)$ in \eqref{beta}.

 An extension of Ramanujan's Master Theorem has been studied by M.~A.~Chaudhry and A.~Qadir \cite{cq}.

\section{Ramanujan's (Earlier) Notebooks}

As we have seen in the preceding sections, many of Ramanujan's theorems and examples on integrals appear in his notebooks \cite{nb}.  Because of their centrality in Ramanujan's vast accomplishments in his theories of theta functions and modular equations, our concentration in this section focuses on elliptic integrals.

The complete elliptic integral of the first kind $K(k)$, $0<k<1$, is defined by
 \begin{equation*}
 K(k):=\int_0^{\pi/2}\df{dt}{\sqrt{1-k^2\sin^2t}}=\df{\pi}{2}\,{_2F_1}\left(\tf12,\tf12;1;k^2\right),
 \end{equation*}
  where the second equality arises from expanding the integrand in a power series and integrating termwise.
 The number $k$ is called the \emph{modulus}.  The function ${_2F_1}$ on the right-hand side is an (ordinary) hypergeometric function, which is defined (more generally) by
\begin{equation}\label{2F1}
{_pF_q}(a_1,a_2,\dots, a_p;b_1,b_2,\dots, b_q;z):=\sum_{n=0}^{\i}\df{(a_1)_n(a_2)_n\cdots(a_p)_n}{(b_1)_n(b_2)_n\cdots(b_q)_n(n!)}z^n,
\end{equation}
 where
$$(a)_0=1, \qquad (a)_n:=a(a+1)(a+2)\cdots(a+n-1),\qquad n\geq 1,$$
and it is assumed that $p$ and  $q$ are chosen so that \eqref{2F1} converges in some domain.
If $\pi/2$ is replaced by another number $v$, $0<v<\pi/2,$ then the integral is called an incomplete elliptic integral.
The integral $K(k)$ is prominent in the theory of the Jacobian elliptic functions sn$(u)$, cn$(u)$, and dn$(u)$, which evidently were not considered by Ramanujan.
 More importantly, for Ramanujan, $K(k)$ plays a central role in his theories of theta functions, class invariants, singular moduli, Eisenstein series and partitions; its importance cannot be overemphasized.     Any statement about an elliptic integral yields a corresponding statement about an ordinary hypergeometric function, and conversely. However, instead of concentrating on hypergeometric functions, the focus here is on the integral themselves, in particular, their transformations and values.

  Elliptic integrals appear at scattered places in Ramanujan's notebooks.
A particularly rich source of identities for elliptic integrals is Section 7 of Chapter 17 in Ramanujan's second notebook \cite{nb}, \cite[pp.~104--117]{III}.

 We begin with the famous \emph{addition theorem} for elliptic integrals.  Let
\begin{equation*}
u:=\int_0^{\alpha}\df{d \varphi}{\sqrt{1-x^2\sin^2\varphi}}, \quad v:=\int_0^{\beta}\df{d \varphi}{\sqrt{1-x^2\sin^2\varphi}}, \quad
w:=\int_0^{\gamma}\df{d \varphi}{\sqrt{1-x^2\sin^2\varphi}}.
\end{equation*}
Ramanujan gave four different conditions for $\alpha$, $\beta$, and $\gamma$ to ensure the validity of the addition theorem \cite[p.~107]{III}
\begin{equation}\label{uvw}
 u+v=w.
 \end{equation}
 In particular, if \cite[Entry 7(viii) (c)]{III},
 \begin{equation}\label{uvw1}
\cot\,\alpha\,\cot\,\beta=\df{\cos\,\gamma}{\sin\,\alpha\,\sin\,\beta}+\sqrt{1-x\sin^2\gamma},
\end{equation}
then \eqref{uvw} holds. The condition \eqref{uvw1} is equivalent to the condition
$$ \df{\text{cn}(u)\,\text{cn}(v)}{\text{sn}(u)\,\text{sn}(v)}=\df{\text{cn}(u+v)}{\text{sn}(u)\,\text{sn}(v)}+\text{dn}(u+v).$$
Although the addition theorem \eqref{uvw} is classical, many of Ramanujan's identities involving elliptic integrals appear to be new.

 In \cite[pp.~108--109]{III}, the two given proofs of the following result are verifications; Ramanujan must have had a more natural proof.
\begin{entry}
If $|x|<1$, then
\begin{equation*}
\df{\pi}{2}\int_0^{\pi/2} \df{d \varphi}{\sqrt{1+x\sin\,\varphi}}=\int_0^{\pi/2}\df{\cos^{-1}(x\,\sin^2\varphi)d\,\varphi}{\sqrt{1-x^2\sin^4\varphi}}.
\end{equation*}
\end{entry}

The following entry is a beautiful theorem, more recondite than the previous theorem. It is a wonderful illustration of Ramanujan's ingenuity and quest for beauty \cite[pp.~111--112]{III}.
\begin{entry}
If $|x|<1$, then
\begin{gather*}
\int_0^{\pi/2}\int_0^{\pi/2}\df{x\sin\,\varphi\,d\theta\,d\varphi}{\sqrt{(1-x^2\sin^2\varphi)(1-x^2\sin^2\theta\,\sin^2\varphi)}}\\
=\df12\left(\int_0^{\pi/2}\df{d\varphi}{\sqrt{1-\tf12(1+x)\sin^2\varphi}}\right)^2-
\df12\left(\int_0^{\pi/2}\df{d\varphi}{\sqrt{1-\tf12(1-x)\sin^2\varphi}}\right)^2.
\end{gather*}
\end{entry}

Despite the fact that Ramanujan's second notebook is a revised edition of the first, there are over 200 claims in the first notebook that cannot be located in the second.  In particular, on page 172 in the first notebook \cite{nb}, two remarkable elliptic integral transformations are recorded \cite[pp.~403--404]{V}. One of them is given in the next entry.
\begin{entry}
Let $0<x<1$, and assume for $0\leq\alpha, \beta\leq \pi/2$ that
$$ \df{1+\sin\,\beta}{1-\sin\,\beta}=\df{1+\sin\,\alpha}{1-\sin\,\alpha}\left(\df{1+x\sin\,\alpha}{1-x\sin\,\alpha}\right)^2.$$
Then,
\begin{equation*}
(1+2x)\int_0^{\alpha}\df{d\theta}{\sqrt{1-x^3\left(\frac{2+x}{1+2x}\right)\sin^2\theta}}=
\int_0^{\beta}\df{d\theta}{\sqrt{1-x\left(\frac{2+x}{1+2x}\right)^3\sin^2\theta}}.
\end{equation*}
\end{entry}

The next two unusual entries are related to elliptic integrals and are found in the unorganized pages of Ramanujan's second notebook \cite[pp.~283, 286]{nb}, \cite[p.~255]{IV}.
\begin{entry}
Let $0\leq\theta\leq\pi/2$ and $0\leq v\leq1$.  Define $\mu$ to be the constant defined by putting $v=1$ and $\theta=\pi/2$ in the definition
\begin{equation}\label{G}
\df{\theta\mu}{2}=\int_0^v\df{dt}{\sqrt{1+t^4}}=:G(v).
\end{equation}
Then,
\begin{equation}\label{el}
2\tan^{-1}v=\theta+\sum_{n=1}^{\i}\df{\sin(2n\theta)}{n\cosh(n\pi)}.
\end{equation}
\end{entry}
Despite its unusual character, \eqref{el} is not too difficult to prove, and follows from the inversion theorem for elliptic integrals.

The integral $G(v)$ has a striking resemblance to the classical lemniscate integral defined next.
 As above, let $0\leq\theta\leq\pi/2$ and $0\leq v\leq1$. Define $\mu$ to be the constant defined by putting $v=1$ and $\theta=\pi/2$ in \eqref{F} below. Then the lemniscate integral $F(v)$ is defined by
\begin{equation}\label{F}
\df{\theta\mu}{\sqrt2}=\int_0^v\df{dt}{\sqrt{1-t^4}}=:F(v)=\sum_{n=0}^{\i}\df{\left(\frac12\right)_nv^{4n+1}}{n!(4n+1)},
\end{equation}
where the right-hand side is a representation for $F(v)$ that arises from expanding the integrand in a binomial series.
Ramanujan offers an inversion formula for the lemniscate integral analogous to \eqref{el}.  Altogether, Ramanujan states ten inversion formulas, six of them for the lemniscate integral \cite[pp.~283, 285, 286]{nb}.  We offer one of them \cite[p.~252]{IV}.  Proofs for all six are given in \cite[245--260]{IV}.
\begin{entry} Let $\theta$ and $v$ be as given in \eqref{F}.  Then,
\begin{gather*}
\log\,v+\df{\pi}{6}-\df12 \log 2+\sum_{n=1}^{\infty}\df{(\tf14)_nv^{4n}}{(\tf34)_n4n}
=\log(\sin\theta)+\df{\theta^2}{2\pi}-2\sum_{n=1}^{\infty}\df{\cos(2n\theta)}{n(e^{2\pi n}-1)}.
\end{gather*}
\end{entry}

 If
$$ v=\df{\sqrt2\,x}{\sqrt{1+x^4}},$$
then
$$ F(v)=\int_0^v\df{dt}{\sqrt{1-t^4}}=\sqrt2\int_0^x\df{dt}{\sqrt{1+t^4}}=\sqrt{2}G(x),$$
which is an important key step in the historically famous problem of  doubling the arc length of the lemniscate.

The lemniscate integral was initially studied by James Bernoulli and Count Giulio Fagn{\'{a}}no.  Raymond Ayoub \cite{ayoub} wrote a very informative article emphasizing its history and importance.  Carl Ludwig Siegel \cite{siegel} considered the lemniscate integral so important that he began his development of the theory of elliptic functions with a thorough discussion of it.

\section{Ramanujan's Lost Notebook}
 On pages 51--53 in his lost notebook \cite{lnb}, Ramanujan states several original, surprising, and unusual integral identities involving elliptic integrals and his theta functions, including
\begin{equation}\label{f}
f(-q):=(q;q)_{\infty}, \qquad |q|<1,
\end{equation}
which, except for a factor of $q^{1/24}$, is Dedekind's eta-function $\eta(\tau)$, where $q=e^{2\pi i\tau}, \tau\in\mathbb{H}.$  Ramanujan's integrals of theta functions are associated with elliptic integrals and modular equations of degrees 5, 10, 14, or 35.  In view of degrees 14 and 35, it is surprising that none of degree 7 are given.  These integral identities  were first proved by S.~Raghavan and S.~S.~Rangachari \cite{rr} using the theory of modular forms, and later by the first author, Heng Huat Chan, and Sen-Shan Huang employing ideas with which Ramanujan would have been familiar \cite{bch}. Proofs for all of Ramanujan's identities can be found in Andrews' and the first author's book \cite[pp.~327--371]{geabcbI}.  Certain proofs depend upon transformations of elliptic integrals found in Ramanujan's second notebook and discussed above.  Differential equations for products or quotients of theta functions are also featured in some proofs.  The first of two examples that we give is associated with modular equations of degree 5.

\begin{entry}\cite[p.~333]{geabcbI}, \cite[p.~52]{lnb} Let $f(-q)$ be defined by \eqref{f}, $\epsilon=(\sqrt5+1)/2$, and
\begin{equation}\label{rrcf}
u:=u(q):=\df{q^{1/5}}{1}\+\df{q}{1}\+\df{q^2}{1}\+\df{q^3}{1}\+\cds, \qquad |q|<1,
\end{equation}
which defines the Rogers--Ramanujan continued fraction.
Then,
\begin{align*}
5^{3/4}\int_0^{q}\df{f^2(-t)f^2(-t^5)}{\sqrt{t}}dt
&=\int_{\cos^{-1}\left((\epsilon u)^{5/2}\right)}^{\pi/2}\df{d\varphi}{\sqrt{1-\epsilon^{-5}5^{-3/2}\sin^2\varphi}}\\
&=\int_0^{2\tan^{-1}\left(5^{3/4}\sqrt{q}f^3(-q^5)/f^3(-q)\right)}\df{d\varphi}{\sqrt{1-\epsilon^{-5}5^{-3/2}\sin^2\varphi}}.
\end{align*}
\end{entry}

To prove the next entry, also associated with modular equations of degree 5, we need a differential equation involving theta functions.

\begin{lemma} Let
$$\lambda:=\lambda(q):=q\df{f^6(-q^5)}{f^6(-q)}.$$
Then
$$q\df{d}{dq}\lambda(q)=\sqrt{q}\,f^2(-q)f^2(-q^5)\sqrt{125\lambda^3+22\lambda^2+\lambda}.$$
\end{lemma}

\begin{entry}\cite[p.~342]{geabcbI}, \cite[p.~52]{lnb} Let $u$ be defined by \eqref{rrcf}.  Then there exists a constant $C$ such that
\begin{equation*}
u^5+u^{-5}=\df{1}{2\sqrt{q}}\df{f^3(-q)}{f^3(-q^5)}\left(C+\int_q^1\df{f^8(-t)}{f^4(-t^5)}\df{dt}{t^{3/2}}+125
\int_0^q\df{f^8(-t^5)}{f^4(-t)}\sqrt{t}\,dt\right).
\end{equation*}
\end{entry}

\noindent(The constant $C$ can be determined, but it is different from that claimed by Ramanujan \cite[pp.~346--347]{geabcbI}.)

The next entry is connected with modular equations of degree 14.

\begin{entry}\cite[pp.~51--52]{lnb}, \cite[p.~359]{geabcbI}
Let
$$v:=v(q):=q\left(\df{f(-q)f(-q^{14})}{f(-q^2)f(-q^7)}\right)^4.$$
Put
$$c=\df{\sqrt{13+16\sqrt2}}{7}.$$
Then
\begin{gather*}
\int_0^q f(-t)f(-t^2)f(-t^7)f(-t^{14})dt=
\df{1}{\sqrt{8\sqrt2}}\int^{\cos^{-1}c}_{\cos^{-1}\left(c\frac{1+v}{1-v}\right)}\df{d\varphi}{\sqrt{1-\frac{16\sqrt2-13}{32\sqrt2}\sin^2\varphi}}.
\end{gather*}
\end{entry}

Our concluding example of Ramanujan's exquisite formulas is an identity linking elliptic integrals and modular equations of degree 35.

\begin{entry}\cite[p.~53]{lnb}, \cite[p.~364]{geabcbI}
If
$$v:=v(q):=q\df{f(-q)f(-q^{35})}{f(-q^5)f(-q^7)},$$
then
\begin{equation*}
\int_0^q t\,f(-t)f(-t^5)f(-t^7)f(-t^{35})dt=
\int_0^v\df{t\,dt}{\sqrt{(1+t-t^2)(1-5t-9t^3-5t^5-t^6)}}.
\end{equation*}
\end{entry}


\begin{thebibliography}{00}



\bibitem{geabcbI}
G.~E.~Andrews and B.~C.~Berndt, \emph{Ramanujan's Lost Notebook},
Part I, Springer, New York, 2005.

\bibitem{geabcbIV}
G.~E.~Andrews and B.~C.~Berndt, \emph{Ramanujan's Lost Notebook},
Part IV, Springer, New York, 2013.

\bibitem{askey}
R.~Askey, \emph{Ramanujan's extensions of the gamma and beta functions}, Amer.~Math.~Monthly \textbf{87} (1980), 346--359.


\bibitem{ayoub}
R.~Ayoub, \emph{The lemniscate and Fagn{\'{a}}no's contributions to elliptic integrals}, Arch.~Hist.~Exact Sci.~\textbf{29} (1984), 131--149.

 \bibitem{quartblms}
 B.~C.~Berndt, \emph{Ramanujan's quarterly reports}, Bull.~London Math.~Soc.~\textbf{16} (1984), 449--489.


\bibitem{I}
B.~C.~Berndt, \emph{Ramanujan's Notebooks}, Part I,
Springer-Verlag, New York, 1985.

\bibitem{II}
B.~C.~Berndt, \emph{Ramanujan's Notebooks}, Part II,
Springer-Verlag, New York, 1989.


\bibitem{rje}
B.~C.~Berndt and R.~J.~Evans, \emph{An integral functional equation of Ramanujan related to the
dilogarithm}, in {\it Number Theory}, Proc.~First Conf.~Canadian Number
Theory Assoc., Banff, R.~A.~Mollin, ed., Walter de Gruyter, Berlin,
1990, pp.~1--5.

\bibitem{III}
B.~C.~Berndt,
\emph{Ramanujan's Notebooks}, Part III,
Springer--Verlag,
New York,
1991.


\bibitem{IV}
B.~C.~Berndt, \emph{Ramanujan's Notebooks}, Part IV,
Springer-Verlag, New York,  1994.

\bibitem{bcbrar}
B.~C.~Berndt and R.~A.~Rankin, \emph{Ramanujan: Letters and Commentary}, Amer.~Math.~Soc., Providence, RI; London Math.~Soc., 1995.

\bibitem{V}
B.~C.~Berndt,
\emph{Ramanujan's Notebooks}, Part V,
Springer--Verlag,
New York, 1998.


\bibitem{lange}
B.~C.~Berndt,  Y.--S.~Choi, and S.--Y.~Kang, \emph{The problems submitted by Ramanujan to the Journal of the Indian
Mathematical Society}, in {\it Continued
Fractions:  From Analytic Number Theory to Constructive Approximation},
B.~C.~Berndt and F.~Gesztesy, eds.,
Contem.~Math.~236, American Mathematical Society, Providence, RI, 1999,
 pp.~15--56.



\bibitem{bch}
B.~C.~Berndt, H.~H.~Chan, and S.--S.~Huang, \emph{Incomplete
elliptic integrals in Ramanujan's lost notebook}, in
\emph{$q$-Series from a Contemporary Perspective}, M.~E.~H.~Ismail
and D.~Stanton, eds., American Mathematical Society, Providence,
RI, 2000, pp.~79--126.


\bibitem{bcbrara}
B.~C.~Berndt and R.~A.~Rankin, \emph{Ramanujan: Essays and Surveys},
American Mathematical Society, Providence, RI, 2001; London
Math.~Soc., 2001.

\bibitem{bcbad}
B.~C.~Berndt and A.~Dixit, \emph{A transformation formula involving the Gamma and Riemann zeta functions in Ramanujan's Lost Notebook}, \emph{The legacy of Alladi Ramakrishnan in the mathematical sciences}, K. Alladi, J. Klauder, C. R. Rao, Eds, Springer, New York, 2010, pp.~199--210.

\bibitem{cq}
M.~A.~Chaudhry and A.~Qadir, \emph{Extension of Hardy's class for Ramanujan's interpolation formula and master theorem with applications}, J.~Inequal.~Appl.~\textbf{52} (2012), 13 pages.

\bibitem{dixit}
A.~Dixit, \emph{Analogues of a transformation formula of Ramanujan}, Int. J. Number Theory~\textbf{7}, No. 5 (2011), 1151-1172.

\bibitem{dixitms}
A.~Dixit, \emph{Modular-type transformations and integrals involving the Riemann $\Xi$-function}, Math.~Student~\textbf{87} Nos. 3-4 (2018), 47--59.

\bibitem{dixzah}
A.~Dixit and A.~Zaharescu, \emph{Ramanujan's paper on Riemann's functions $\xi(s)$ and $\Xi(t)$ and a transformation from the Lost Notebook}, to appear in the \emph{Encyclopedia of Srinivasa Ramanujan and His Mathematics}.

\bibitem{gasper}
G.~Gasper and M.~Rahman, \emph{Basic Hypergeometric Series}, 2nd ed., Cambridge University Press, Cambridge, 2004.



\bibitem{hardy1}
G.~H.~Hardy, \emph{Note on the function $\int_x^{\infty}e^{\frac12 (x^2-t^2)}dt$}, Quart.~J.~Math.~\textbf{35} (1904), 193--207.

\bibitem{ghh}
G.~H.~Hardy, \emph{Note by Mr. G.H.~Hardy on the preceding paper}, Quart.~J.~Math.~\textbf{46} (1915),
260--261.

\bibitem{hardyw}
G.~H.~Hardy, \emph{Mr. S.~Ramanujan's mathematical work in England}, J. Indian Math. Soc.~\textbf{9} (1917), 30--45; (also a Report to the University of Madras, 1916, privately printed).

 \bibitem{ram}
 G.~H.~Hardy, \emph{Ramanujan}, Chelsea, New York, 1978.



\bibitem{mustafy}
A.~K.~Mustafy, \emph{A new representation of Riemann's zeta function and some of its consequences}, Kongelige Norske Videnskabers Selskabs Forhandlinger \textbf{39} (1966), 96--100.

\bibitem{preece}
C.~T.~Preece, \emph{Theorems stated by Ramanujan (X)}, J.~London Math.~Soc.~\textbf{6} (1931), 22--32.

\bibitem{rr}
S.~Raghavan and S.~S.~Rangachari,
\emph{On Ramanujan's elliptic integrals and modular identities} in
\emph{Number Theory and Related Topics}, Oxford University Press, Bombay,
1989, pp.~119--149.

\bibitem{295}
S.~Ramanujan, \emph{Question 295}, J.~Indian Math.~Soc.~\textbf{3} (1911), 128.


\bibitem{463}
S.~Ramanujan, \emph{Question 463}, J.~Indian Math.~Soc.~\textbf{5} (1913), 120.



\bibitem{7}
S.~Ramanujan, \emph{On the integral
${\int_0^x\tf{\tan^{-1}t}{t}dt}$}, J.~Indian
Math.~Soc.~\textbf{7} (1915), 93--96.

\bibitem{11}
S.~Ramanujan, \emph{Some definite integrals},
Mess.~Math.~\textbf{44} (1915), 10--18.

\bibitem{12}
S.~Ramanujan, \emph{Some definite integrals connected with Gauss's
sums}, Mess. Math.~\textbf{44} (1915), 75--85.

\bibitem{riemann}
S.~Ramanujan, \emph{New expressions for Riemann's functions
$\xi(s)$ and $\Xi(t)$}, Quart.~J.~Math.~\textbf{46} (1915),
253--260.

\bibitem{22}
S.~Ramanujan, \emph{Some definite integrals}, Proc.~London Math.~Soc.~\textbf{17} (1918), Records for 17 January 1918,

\bibitem{783}
S.~Ramanujan, \emph{Question 783}, J.~Indian Math.~Soc.~\textbf{10} (1918), 397--399.

\bibitem{23}
S.~Ramanujan, \emph{Some definite integrals}, J.~Indian
Math.~Soc.~\textbf{11} (1919), 81--87.

\bibitem{27}
S.~Ramanujan, \emph{A class of definite integrals},
Quart.~J.~Math.~\textbf{48} (1920), 294--310.





\bibitem{cp}
S.~Ramanujan, \emph{Collected Papers}, G.~H.~Hardy, P.~V.~Seshu Aiyar, and B.~M.~Wilson, eds., Cambridge University Press,
Cambridge, 1927; reprinted by Chelsea, New York, 1962; reprinted
by the American Mathematical Society, Providence, RI, 2000.


\bibitem{nb}
S.~Ramanujan, \emph{Notebooks of Srinivasa Ramanujan} (2 volumes), Tata Institute of
Fundamental Research, Bombay, 1957; second ed., 2012.


\bibitem{lnb}
S.~Ramanujan, \emph{The Lost Notebook and Other Unpublished
Papers}, Narosa, New Delhi, 1988.


\bibitem{siegel}
C.~L.~Siegel, \emph{Topics in Complex Function Theory}, Vol.~1, Wiley, New York, 1969.

\bibitem{titch}
E.~C.~Titchmarsh, \emph{The Theory of the Riemann Zeta Function}, Clarendon Press, Oxford, 1986.

\bibitem{watson}
G.~N.~Watson, \emph{Theory of Bessel Functions}, Cambridge University Press, second ed., London, 1966.

\end{thebibliography}
\end{document}